\documentclass[a4paper, 11pt]{article}
\pdfoutput=1
\usepackage[margin=1in]{geometry}
\usepackage{pslatex}
\usepackage{fancyhdr}
\usepackage{graphicx}
\usepackage{amsmath}
\usepackage{amsfonts}
\usepackage{amssymb}
\usepackage{float}
\usepackage{tikz}
\usepackage{graphicx}
\usepackage{amsmath}
\usepackage{amsfonts}
\usepackage{amssymb}
\usepackage{parskip}
\renewcommand\abstract{%
  \textbf{Abstract}
  \begingroup
  \let\textit\relax
}
\endabstract{\vspace{1em}\par}
\makeatother
\title{\Large \textbf{A NUMERICAL MULTISCALE METHOD FOR FIBER NETWORKS}}
\author{Morgan G\"{o}rtz$^1$ \hspace{0.4 cm} Gustav Kettil$^1$ \hspace{0.4 cm} Axel M\aa lqvist$^2$ \hspace{0.4 cm} Andreas Mark$^1$ \hspace{0.4 cm} Fredrik Edelvik$^1$}
\begin{document}
\maketitle
\abstract{
\noindent Fiber network modeling can be used for studying mechanical properties of paper \cite{Kulachenko}. The individual fibers and the bonds in-between constitute a detailed representation of the material. However, detailed microscale fiber network models must be resolved with efficient numerical methods. In this work, a numerical multiscale method for discrete network models is proposed that is based on the localized orthogonal decomposition method \cite{LOD}. The method is ideal for these network problems, because it reduces the maximum size of the problem, it is suitable for parallelization, and it can effectively solve fracture propagation. 

The problem analyzed in this work is the nodal displacement of a fiber network given an applied load. This problem is formulated as a linear system that is solved by using the aforementioned multiscale method. To solve the linear system, the multiscale method constructs a low-dimensional solution space with good approximation properties \cite{DLOD,Phd}. The method is observed to work well for unstructured fiber networks, with optimal rates of convergence obtainable for highly localized configurations of the method.
}
\section{INTRODUCTION}
Network models are used in multiple fields of science. Examples of problems modeled using network models are traffic networks, vascular networks, pours media and fibrous materials. In these models, the micro-scale properties define macro-scale behavior. To handle the enormous computational complexity of huge network models some form of upscaling, homogenization, is normally performed. Rudimentary averaging techniques have shown troublesome, leading to convergence plateaus and erroneous solutions. Therefore numerical multiscale methods have been developed, during the last decades, to provide reliable solutions to these models. 

\begin{figure}[H]
\centering
\includegraphics[width=0.3\linewidth]{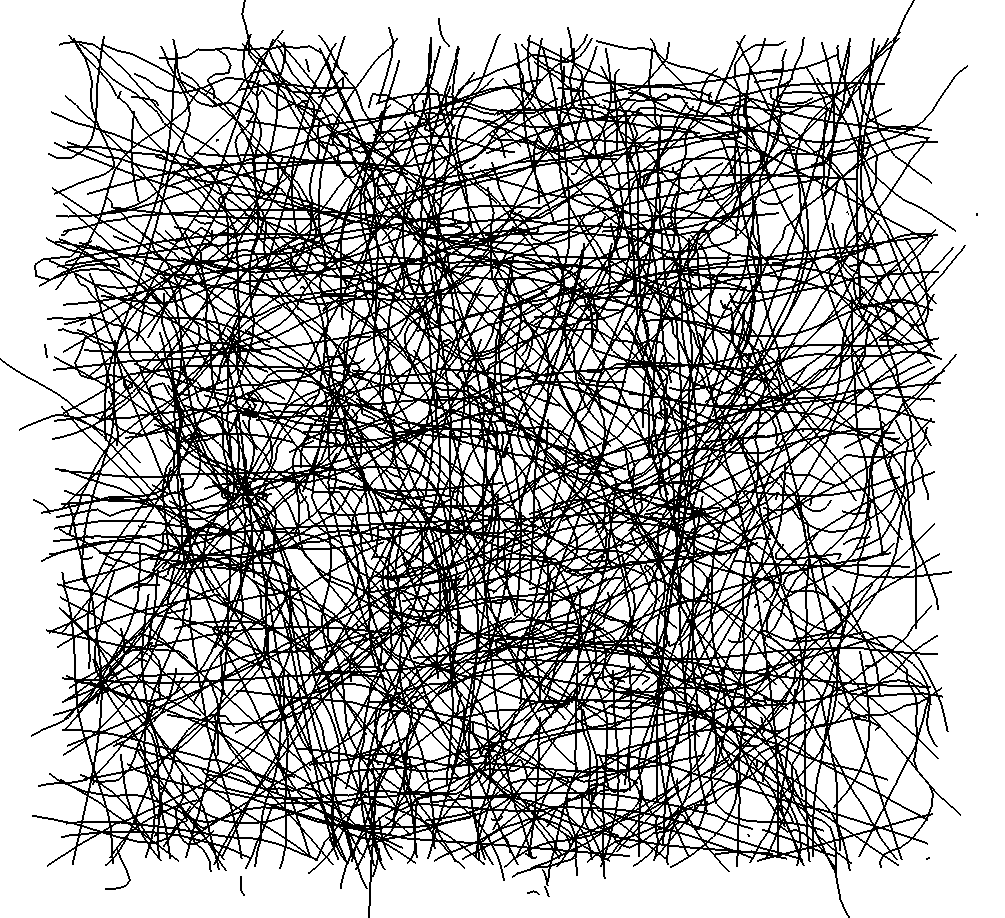}
\caption{Fiber network generated in a paper manufacture simulation. In the simulation, fiber suspension flows onto a woven mesh where the paper is formed.}
\label{fig:simulation}
\end{figure}

There is a vast amount of literature on multiscale methods for solving partial differential equations with rapidly varying coefficients. Two examples are the Heterogeneous Multiscale Method (HMM) \cite{HMM} and the Multiscale Finite Element Method (MsFEM) \cite{MsFEM}. More recently the Localized Orthogonal Decomposition (LOD) \cite{LOD} has gained popularity since it is designed to handle non-periodic and non scale separated data. In \textit{Numerical upscaling of discrete network models} \cite{DLOD} LOD is first applied to a fiber network model. This is the approach used in this paper. In particular, we will apply it to micro-scale paper models with realistic geometries in two dimensions. There are several other works on multiscale methods for network models including Ewing et al. on heat connectivity of network materials \cite{HeatNetwork}, Chu et al. \cite{porous} on porous media, and Tomas et al \cite{WetPaper} on edge wicking in multi-ply paperboard. \\
\hrule
$^1$ Computational Engineering and Design, Fraunhofer-Chalmers Centre, Chalmers Science Park, 412 88 Gothenburg, Sweden. \\
$^2$ Department of Mathematical Sciences, Chalmers University of Technology and University of Gothenburg, 412
96 Gothenburg, Sweden.
\newpage
A sheet of paper can be modeled as a network structure. The models analyzed in this work are inspired by the microscopic structure of a sheet of paper. Paper is constructed by a fiber suspension flowing onto a moving woven mesh. The individual fibers in the suspension bond together on the weave to form a continuous sheet of paper. This process has been simulated \cite{paperlaying} and a cropped generated network from a simulation is presented in Figure \ref{fig:simulation}. In future work, these three-dimensional network geometries will be used to simulate the paper's structural properties. Not all paper fibers have the same structural properties, not all bonds have the same holding strength, so to represent these variations individual fiber strengths are randomized. In this article, it is observed that stiff randomized fiber coefficients and fiber-based network geometries can be homogenized optimally in two dimensions using the LOD method.
\section{PROBLEM FORMULATION} \label{sec:problem}
\begin{figure}[H]
\centering
    \includegraphics[width=.32\textwidth]{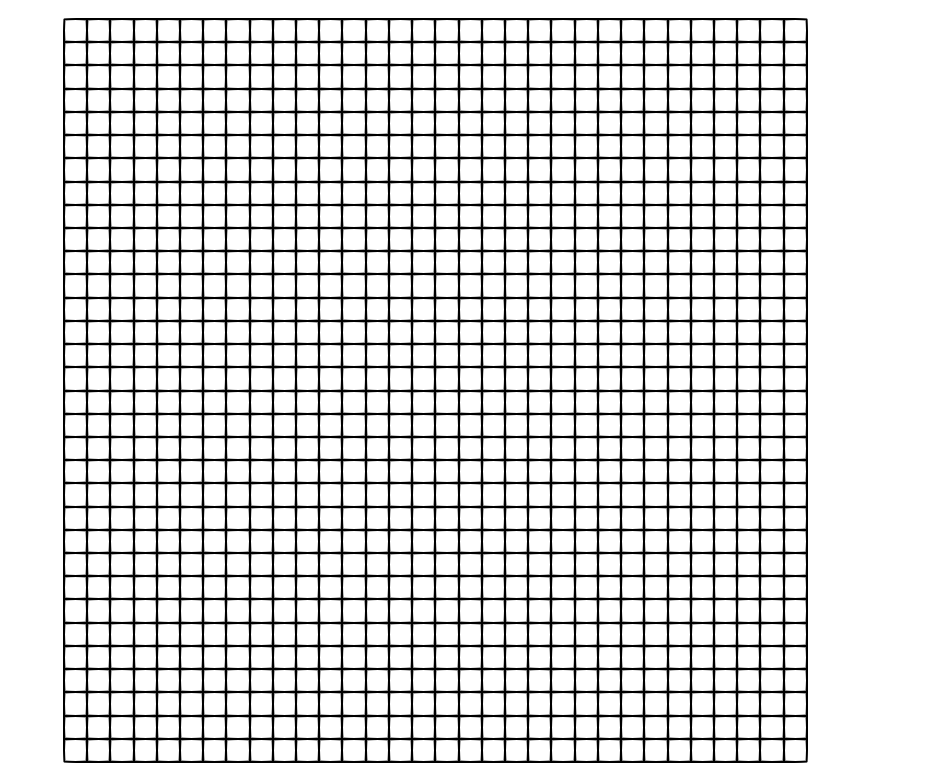}
    \includegraphics[width=.32\linewidth]{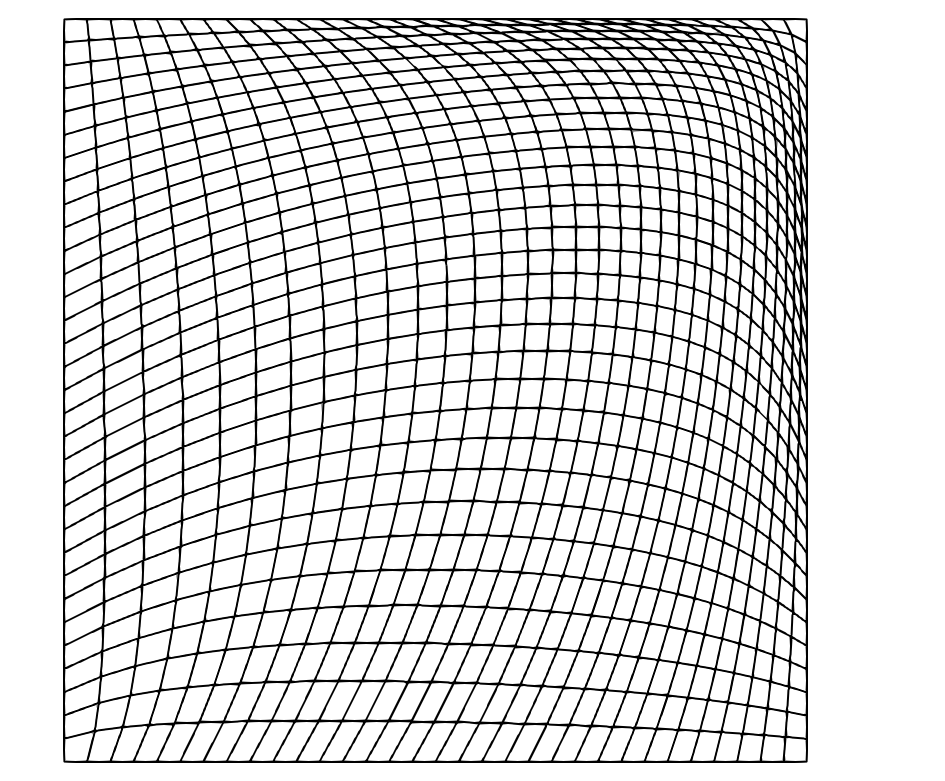}
    \includegraphics[width=.32\linewidth]{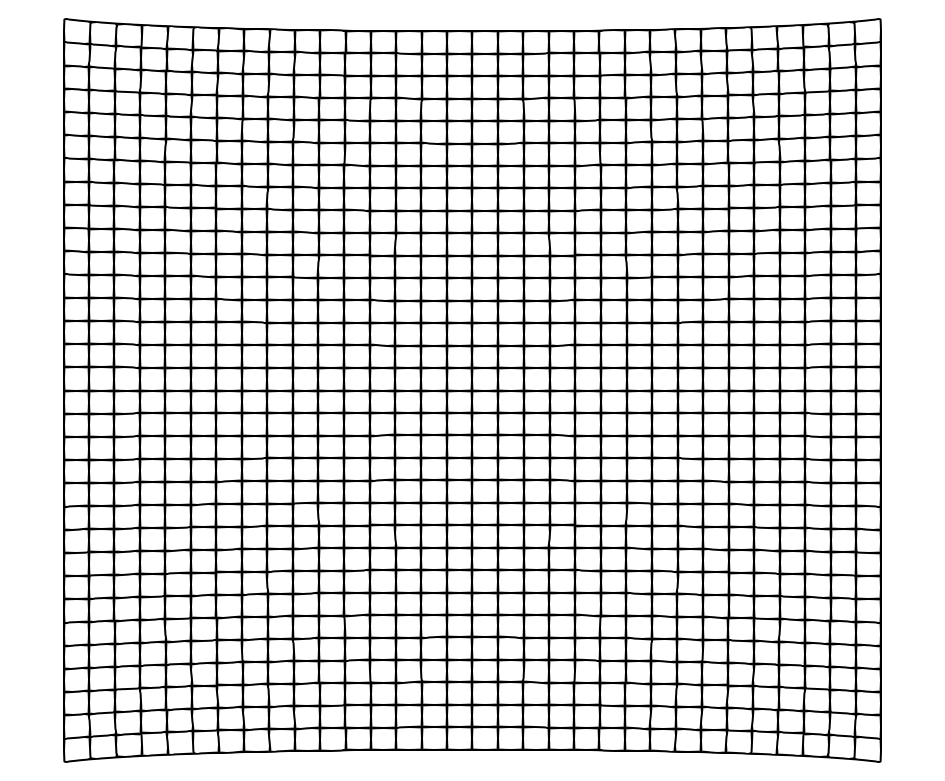}
\caption[2]{Visualization of the network displacement solutions for the two different problems. The left network is the original network. In the middle, the solution when a constant upward diagonal force is applied. The right network is the solution to the problem with a fixed left boundary and a constant displacement on the right boundary.}
\label{fig:diagonal}
\end{figure}

Two different displacement problems are analyzed. The first problem is to find the displacement of a square network model with a fixed boundary given an applied force. In this problem, the force is a constant diagonally upward oriented force scaled with respect to the stiffness of the model. In the second problem, the displacement of a square network model with a fixed left boundary and a 10\% displacement of the right boundary is solved. No force is applied in the second problem. Visualizations of example displacement solutions for both of the problems can be found in Figure \ref{fig:diagonal}.
\subsection{Mathematical formulation} 
The problems are formulated as linear systems on the form:
\begin{equation} \label{eq:problem}
    \begin{split}
        Ku &= F, \\
        u = [u_{1},u_{2}, \hdots, u_{n}]^T, &\ F= [F_{1},F_{2}, \hdots,F_{n}]^T,
    \end{split}
\end{equation}
where $n$ is the number of network nodes, $K\in\mathbb{R}^{2n\times2n}$ connectivity, $u_{i}=[u_{i,x},u_{i,y}]$ are the two-directional nodal displacements, and $F_{i}=[F_{i,x},F_{i,y}]$ are the directional applied forces in the $i$:th network node. The connectivity matrix, $K$, is composed of three linear force-displacement relations. 

The first is edge extension, defining the forces arising in the change of edge length. Given an edge, $e$, the edge extension force in one of the edges nodes, $i$, is defined by:
\begin{equation} \label{eq:edge_extension}
\text{FE}_{i,\text{e}} = k_{\text{e}}\frac{a_{\text{e}}}{L_{\text{e}}}\Delta L_{\text{e}}^id_{\text{e}}^i,
\ \Delta L_{\text{e}}^i = (u_i^{\text{other}}-u_i)d_{\text{e}}^i,
\end{equation} 
where $k_{\text{e}}$, $a_{\text{e}}$, $L_{\text{e}}$, and $d_{\text{e}}^i$ are extension stiffness, cross section area, edge length, and the $i$:th directional normalized edge direction respectively. The $u_i^{\text{other}}$ term is the displacement of the other edge node. 

\begin{figure}[H]
\centering
    \begin{tikzpicture}[scale = 1.3]
    \draw[dashed,-] (0,0) -- (2,2) node[anchor=south]{};
    \draw[thick,-] (0.5,0.5) -- (1.5,1.5) node[anchor=south]{};
    \draw[thick,<->] (0.25,0.25) -- (1.75,1.75) node[anchor=south]{};
    \filldraw (0.5,0.5) circle (1pt);
    \filldraw (1.5,1.5) circle (1pt);
    \filldraw (0,0) circle (1pt);
    \filldraw (2,2) circle (1pt);
    \end{tikzpicture}
    \begin{tikzpicture}[scale = 1.3]
    \draw[dashed,-] (1,0) -- (2,2) node[anchor=south]{};
    \draw[dashed,-] (1,0) -- (0,2) node[anchor=south]{};
    \draw[thick,-] (1,0) -- (2.5,1.7) node[anchor=south]{};
    \draw[thick,-] (1,0) -- (-0.5,1.7) node[anchor=south]{};
    \draw[thick,->] (-0.5,1.7) -- (-0.2,2) node[anchor=south]{};
    \draw[thick,->] (2.5,1.7) -- (2.2,2) node[anchor=south]{};
    \filldraw (1,0) circle (1pt);
    \filldraw (2,2) circle (1pt);
    \filldraw (0,2) circle (1pt);
    \filldraw (-0.5,1.7) circle (1pt);
    \filldraw (2.5,1.7) circle (1pt);
    \draw[thick,->] (1,0) -- (0.7,-0.3) node[anchor=south]{};
    \draw[thick,->] (1,0) -- (1.3,-0.3) node[anchor=south]{};
    \end{tikzpicture}
    \begin{tikzpicture}[scale = 1.3]
    \draw[dashed,-] (0,1.3) -- (0,2) node[anchor=south]{};
    \draw[thick,-] (0,0) -- (2,0) node[anchor=south]{};
    \draw[thick,-] (0,0) -- (0,1.3) node[anchor=south]{};
    
    \draw[thick,->] (0,1.3) -- (0,1.6) node[anchor=south]{};
    \draw[thick,->] (2,0) -- (2.3,0) node[anchor=south]{};
    \draw[thick,->] (0,0) -- (0,-.3) node[anchor=south]{};
    \draw[thick,->] (0,0) -- (-.3,0) node[anchor=south]{};
    \filldraw (2,0) circle (1pt);
    \filldraw (0,2) circle (1pt);
    \filldraw (0,0) circle (1pt);
    \filldraw (0,1.3) circle (1pt);
    \end{tikzpicture}
    \caption{The three force-displacement relations defining the network connectivity matrix. Solid lines are displaced fibers, dashed lines are original fiber positions, nodes are network nodes and arrows are forces. The left is edge extension, the middle is angular deviation, and the right is Poisson effect. Edge extension acts on edges and gives them a spring-like property. Angular deviation and Poisson effect act on edge pairs and add joint stiffness.}
    \label{fig:forces}
\end{figure}
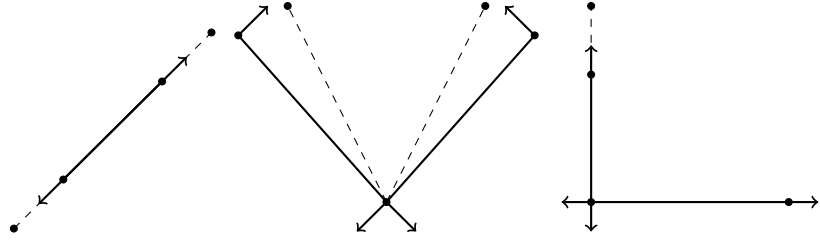

The second force-displacement relation is angular deviation. The definition of angular deviation is as follows: Given an edge pair (ep) composed of two edges, $e_i$ between nodes $(i,j)$ and $e_j$ between nodes $(j,l)$, then the angular deviation force of the edge pair in the node $k$ is defined by:
\begin{equation}\label{eq:angular_devation}
\text{FA}_{k,\text{ep}} = \begin{cases}
-\kappa_{\text{ep}}\dfrac{V_{\text{ep}}}{L_{e_k}}\Delta\theta_{\text{ep}} n^j_{e_k}, \ k = i,l\\
-\text{FA}_{i,\text{ep}} -\text{FA}_{l,\text{ep}}, \ k = j
\end{cases}
\end{equation}
where $\kappa_{\text{ep}}$ is the angular stiffness of the edge pair, $V_{\text{ep}}$ is the connection volume, $\Delta\theta_{\text{ep}}$ is the angular change, and $n^j_{e_k}$ is a directional outward normal of the edge $e_k$. The angular change, $\Delta\theta_{\text{ep}}$, is linearized under the assumption that it is small and is approximated by:
\begin{equation}\label{eq:linang}
    \Delta \theta_{\text{ep}}= \delta\theta_{e_i} + \delta\theta_{e_l}, \ \delta\theta_{e_k} = \frac{(u_k-u_j)}{L_{e_k}} n^j_{e_k}.
\end{equation}
The third and final force-displacement relation is the Poisson effect. It takes a similar form as angular deviation as it also acts on edge pairs. The Poisson effect is defined as follows: Given an edge pair (ep) composed of two edges, $e_i$ between nodes $(i,j)$ and $e_j$ between nodes  $(j,l)$, then the Poisson effect of the edge pair in the node $k$ is defined by:
\begin{equation} \label{eq:poison_effect}
\text{FP}_{k,\text{ep}} = \begin{cases}
-\eta_{\text{ep}}\dfrac{a_{e_k}}{L_{e_k}} \left(\Delta L_{e_k}^j+\gamma_{ep} w_{e_k^\text{other}} \dfrac{\Delta L_{e_k^\text{other}}}{2L_{e_k^{\text{other}}}}|n_{e_k}^j\cdot d_{e_k^{\text{other}}}^i|\right) d^j_{e_k}, \ k = i,l\\
-\text{FP}_{i,\text{ep}} -\text{FP}_{l,\text{ep}}, \ k = j
\end{cases},
\end{equation}
where $\eta_{ep}$ and $\gamma_{ep}$ are Poisson coefficients, $w_{e}$ is the width of the edge $e$, and $e^{\text{other}}_k$ means the edge not containing the node $k$. Visualizations of these force-displacement relations are presented in Figure \ref{fig:forces}. The connectivity matrix, $K$, in \eqref{eq:problem} is then constructed by combining the three force-displacement relations \eqref{eq:edge_extension}, \eqref{eq:angular_devation} linearized with \eqref{eq:linang}, and \eqref{eq:poison_effect} for all edges and edge pairs into a linear system. 
\section{NETWORK MODEL SPECIFICATION} \label{sec:model}
The fiber network model is broken down into two parts: network type and coefficient type. The network type of a fiber network defines how nodes, edges, and edge pairs are placed in the network and the network's coefficient type defines what coefficients are used in \eqref{eq:edge_extension}, \eqref{eq:angular_devation}, and \eqref{eq:poison_effect} for each edge and edge pair. In this work three network types and three coefficient types are analyzed. 

\begin{figure}[H]
    \centering
    \includegraphics[width=.30\linewidth]{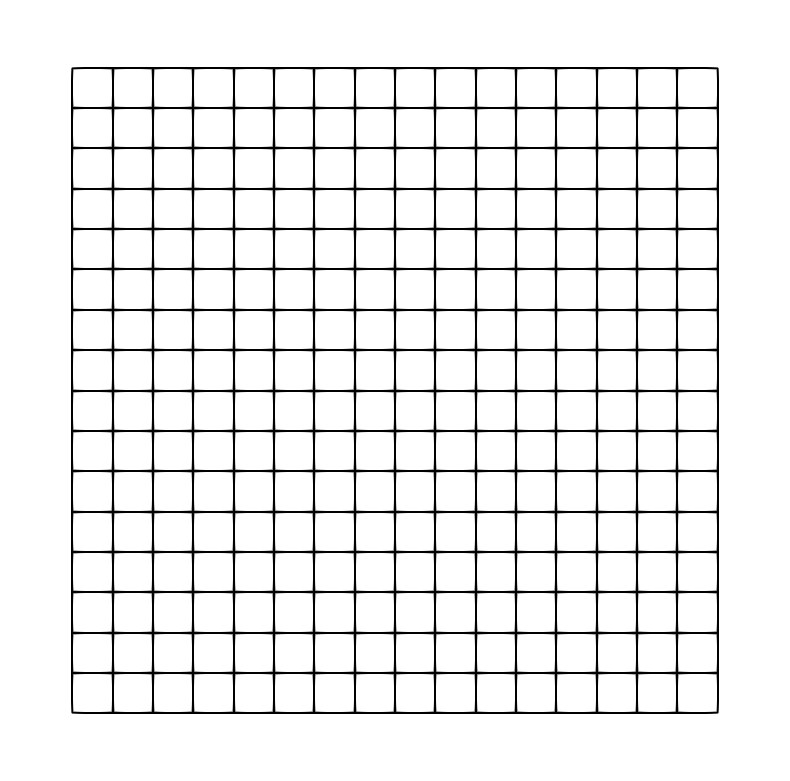}
    \includegraphics[width=.30\linewidth]{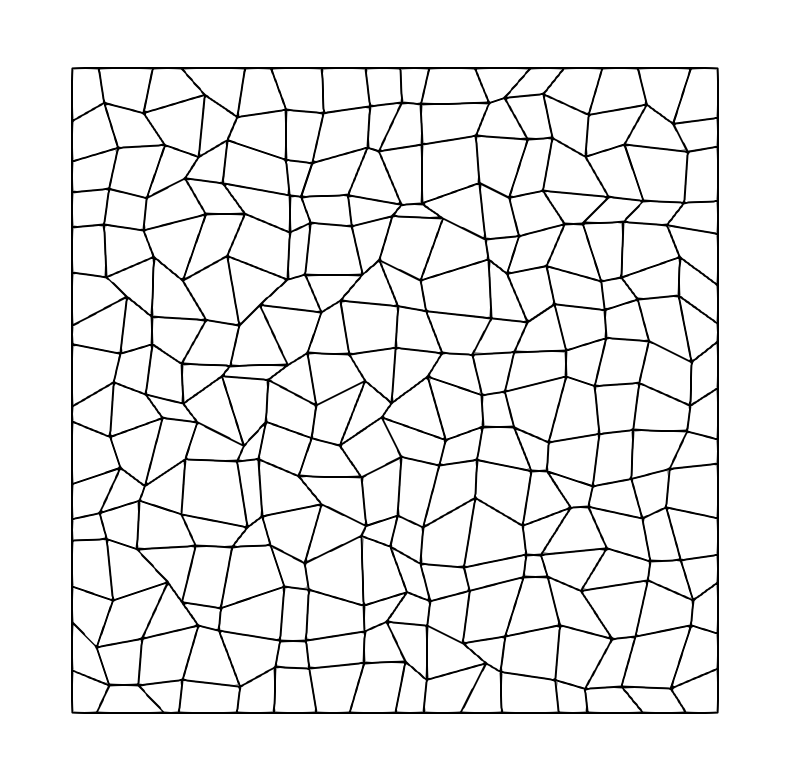}
    \includegraphics[width=.30\linewidth]{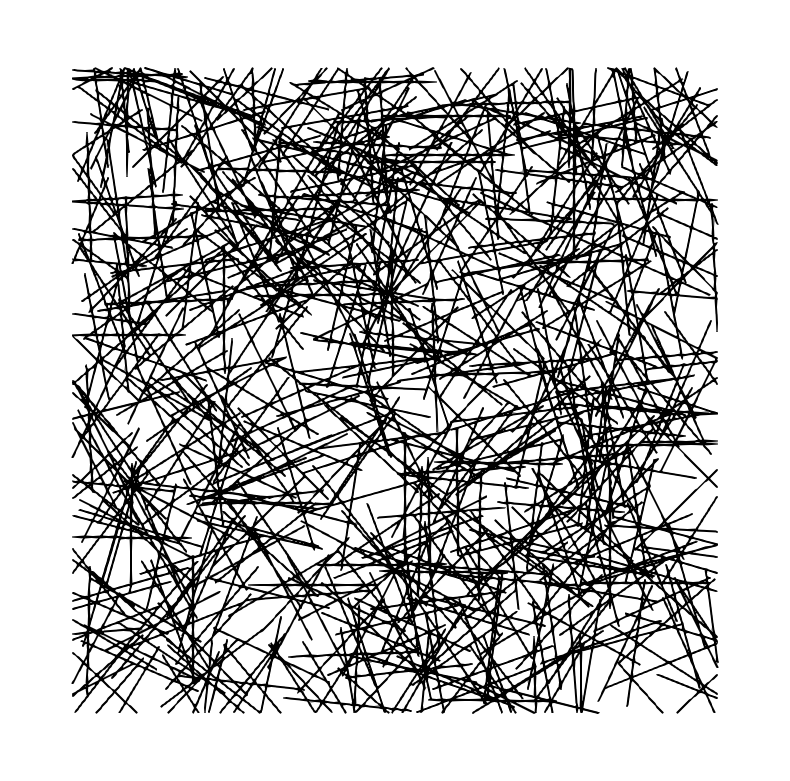}
    \caption{Network types analyzed. From left to right there is a structured, a perturbed, and an unordered fiber network.}
    \label{fig:networks}
\end{figure}

The three network types are structured networks, perturbed networks, and unordered fiber networks. A structured network is an equidistant grid and a perturbed network is a structural network with its nodes perturbed. The unordered fiber network is generated differently. First, fibers are defined as one-dimensional objects with set lengths subdivided into a set of edges. Then a quantity of these fibers is placed at random with their intersections connected with nodes and edge pairs. This process is meant to simulate the paper forming process and has a similar structure to the network presented in Figure \ref{fig:simulation}. Examples of these three network types can be found in Figure \ref{fig:networks}.
For each component, edge or edge pair, force coefficients have to be chosen. Three strategies for choosing the force coefficients are analyzed. The first is homogeneous coefficients. It has been proved \cite{DLOD} that the problem that arises from using a structured network and a specific choice of homogeneous coefficients is equivalent to a finite difference scheme of the linear elasticity equation. The second scheme is randomized coefficients. In this scheme, each edge and edge pair gets different coefficients based on a uniform distribution. This is meant to represent the randomized structural properties of the individual fibers and their adhesion to each other. The last strategy analyzed is fiber coefficients. In this scheme, the edges are meant to have properties similar to real paper fibers, where both uniform and randomized fiber coefficients are analyzed. This scheme creates rigid networks. A structured or perturbed network with fiber coefficients is equivalent to stiff homogeneous or stiff random coefficients. For unordered fiber networks, there are different edge pair coefficients for edge pairs defining the structural properties of fibers and edge pairs defining the structural properties of adhesion between fibers.
\section{METHOD}\label{sec:method}
The Localized Orthogonal Decomposition (LOD) method solves symmetric and positive definite linear systems such as \eqref{eq:problem}. First, the system is reformulated on its weak form:
\begin{equation} \label{eq:weak}
   \text{Find } u \in V: \ v^TKu = v^T F, \  \forall v \in V, 
\end{equation}  
where $V\subset \mathbb{R}^{2n}$ is constructed by removing any vector from $\mathbb{R}^{2n}$ that represents a solution breaking a fixed boundary condition. The LOD method finds an approximation, $u_{ms}$, in a low dimensional multiscale space, $V_{ms}$, by solving:
\begin{equation} \label{eq:weakms}
   \text{Find } u_{ms} \in V_{ms}: \ v^TKu_{ms} = v^T F, \  \forall v \in V_{ms}.
\end{equation} 
This multiscale space is constructed by first covering the network model with a grid. In this paper, two-dimensional equidistant square grids are used. On this square grid, shape functions are defined in each coarse grid node similar to the finite element method. These shape functions are then evaluated in each network node and the span of these discrete shape functions, $\lambda_i$, defines a coarse space $V_H$. Here the bilinear basis functions are chosen as shape functions. There are several ways to define a fine space, $W$, but in this work the fine space is defined as the $l_2$ orthogonal complement of $V_H$ in $V$. This fine space paired with $V_H$ forms a splitting of $V$, $V=V_H\oplus W$, because $(V,(\cdot,\cdot)_{l_2})$ is a vector space. The multiscale space is then defined by the $K$ orthogonal complement of $W$:
$$V_{ms} = \{v\in V: \ w^TKv = 0, \forall w\in W \}$$
and because $K$ is symmetric and coercive (positive definite) $V_{ms}$ and $W$ also forms a splitting of $V$. 

\begin{figure}[H]
    \centering
    \includegraphics[clip, trim= 2.3cm 1.5cm 2.3cm 4cm, width=0.7\textwidth]{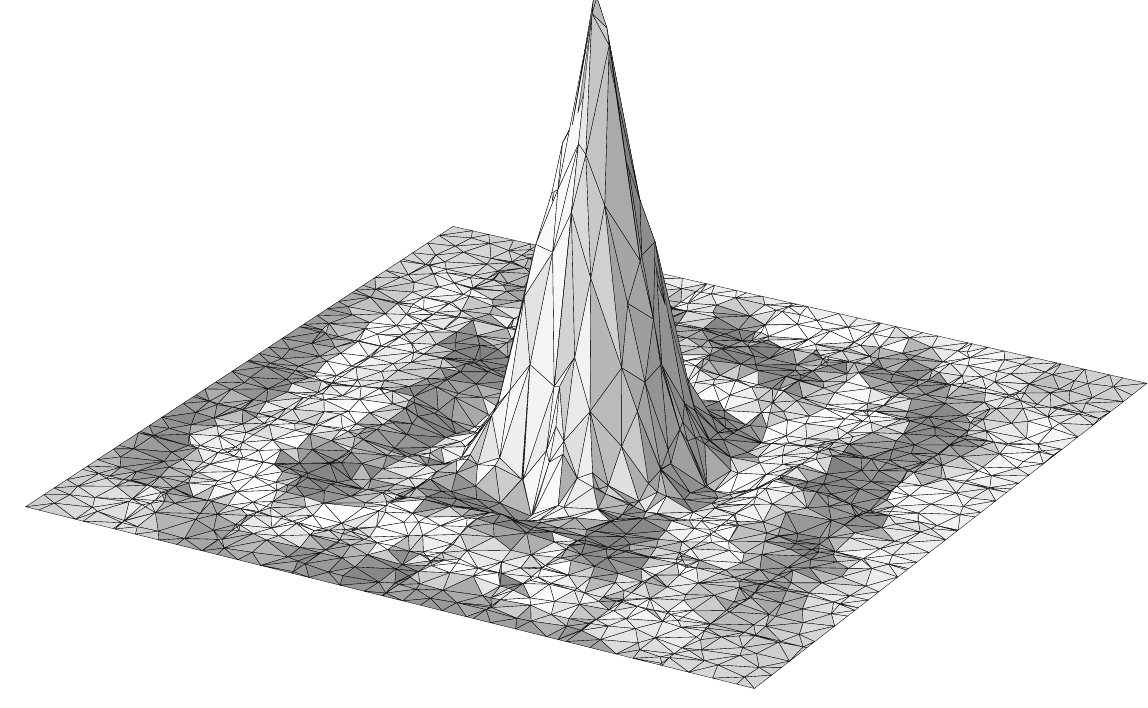}
    \caption{A modified bilinear shape function. The modified shape functions form a basis for the multiscale space and have exponential decay, an essential property to be able to find the modifications locally.}
    \label{fig:modification}
\end{figure}

The method finds a solution by modifying each shape function to construct a basis for the multiscale space $V_{ms}$. To find these modifications, $\phi_i\in W$, the following equations are solved:
$$w^TK(\lambda_i-\phi_i) = 0 \Leftrightarrow w^TK\lambda_i = w^TK\phi_i, \ \forall w\in W.$$
These modifications do not have compact support as they live on the entire domain, however they have exponential decay. This decay allows them to be found locally by assuming they are zero outside a given radius. In Figure \ref{fig:modification}, a modified bilinear shape function, $\lambda_i-\phi_i$, is presented. These modified shape functions form a basis for $V_{ms}$ and are used to solve \eqref{eq:weakms} to get the multiscale approximation. 
\section{NUMERICAL RESULTS}
The multiscale method presented in Section \ref{sec:method} is used to solve the problems presented in Section \ref{sec:problem} on the models presented in Section \ref{sec:model}. Several optimization techniques \cite{APP} are used in the implementation of the LOD method and the method is written in C++. In the experiments the exact solution, $u$, is found by solving the linear system \eqref{eq:problem} and is used to evaluate the error of the multiscale approximation $u_{ms}$. The two errors analyzed for each experiment are:
$$\|u-u_{ms}\|_{l_2} \text{ and } \ \|u-u_{ms}\|_E, $$
where  $\|v\|_E = \sqrt{v^TKv}$ is the energy norm. For the method to be optimal two constants, $C_1$ and $C_2$, should exists such that:
\begin{equation} \label{eq:ineq}
\|u-u_{ms}\|_{l_2} \leq C_1 H^2 \text{ and } \ \|u-u_{ms}\|_E\leq C_2 H,
\end{equation}
where $H$ is the largest diameter of the elements in the coarse grid. For comparison the errors presented in this section are normalized with their corresponding exact solutions. 
Three different setups are analyzed in this work. The first setup is a perturbed network type with uniform and random coefficients, the second is a structured and perturbed network type with stiff homogeneous coefficients, and the third is an unordered fiber network with both homogeneous and random fiber coefficients. All of the setups have some common properties. In each setup, the domain is a square and on this square an equidistant square coarse grid is applied given an element width $H$. This coarse grid is composed of $(m+1) \times (m+1)$ nodal points and the corresponding shape functions on the coarse grid are chosen to be the bilinear basis functions. To make the LOD method localized, the modifications vanish $1.5H\log(m)$ away from each respective nodal point. The inequalities \eqref{eq:ineq} are numerically confirmed in all of these setups.
\subsection{Perturbed network}
\begin{figure}[H]
    \centering
    \includegraphics[width=.43\linewidth]{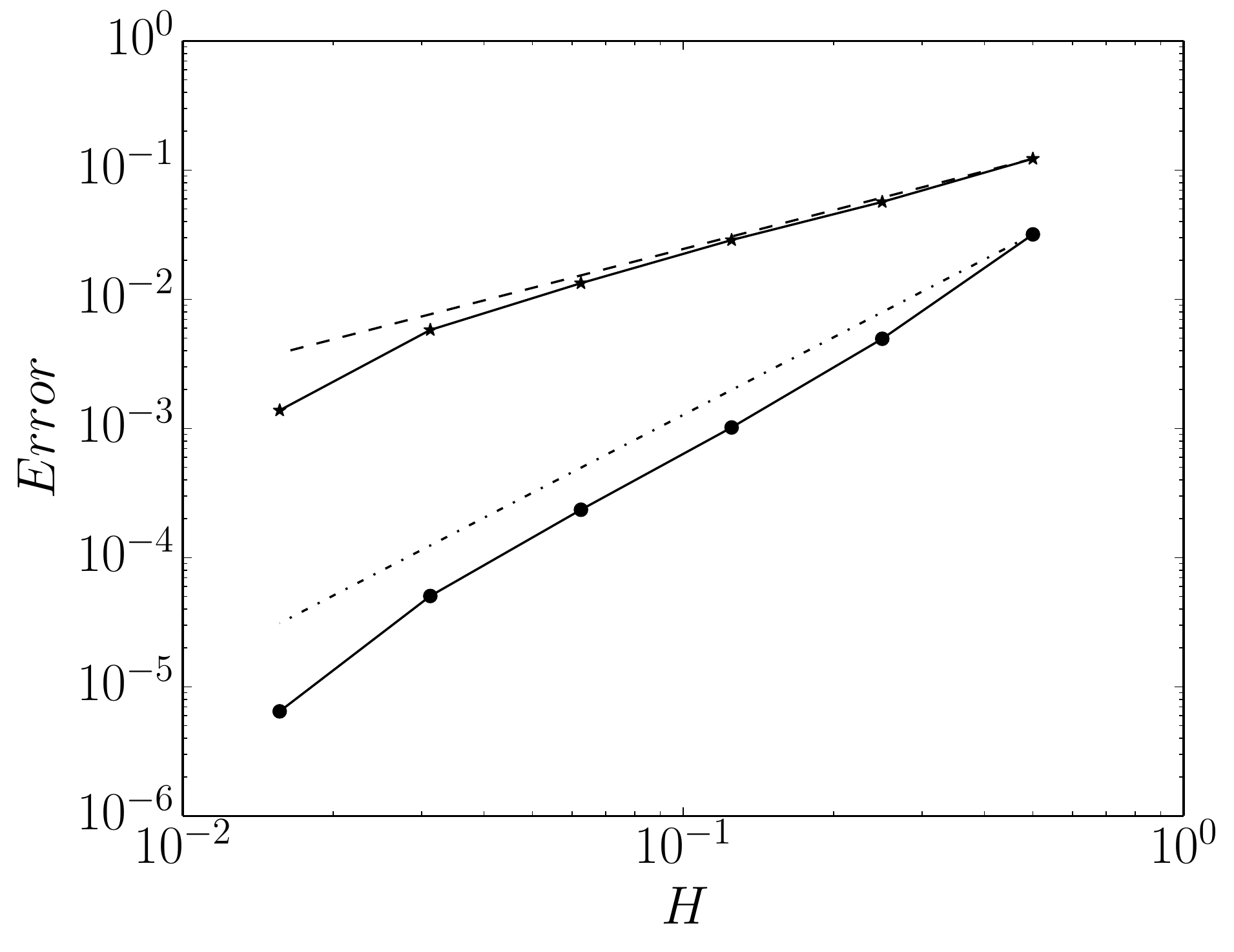}
    \includegraphics[width=.43\linewidth]{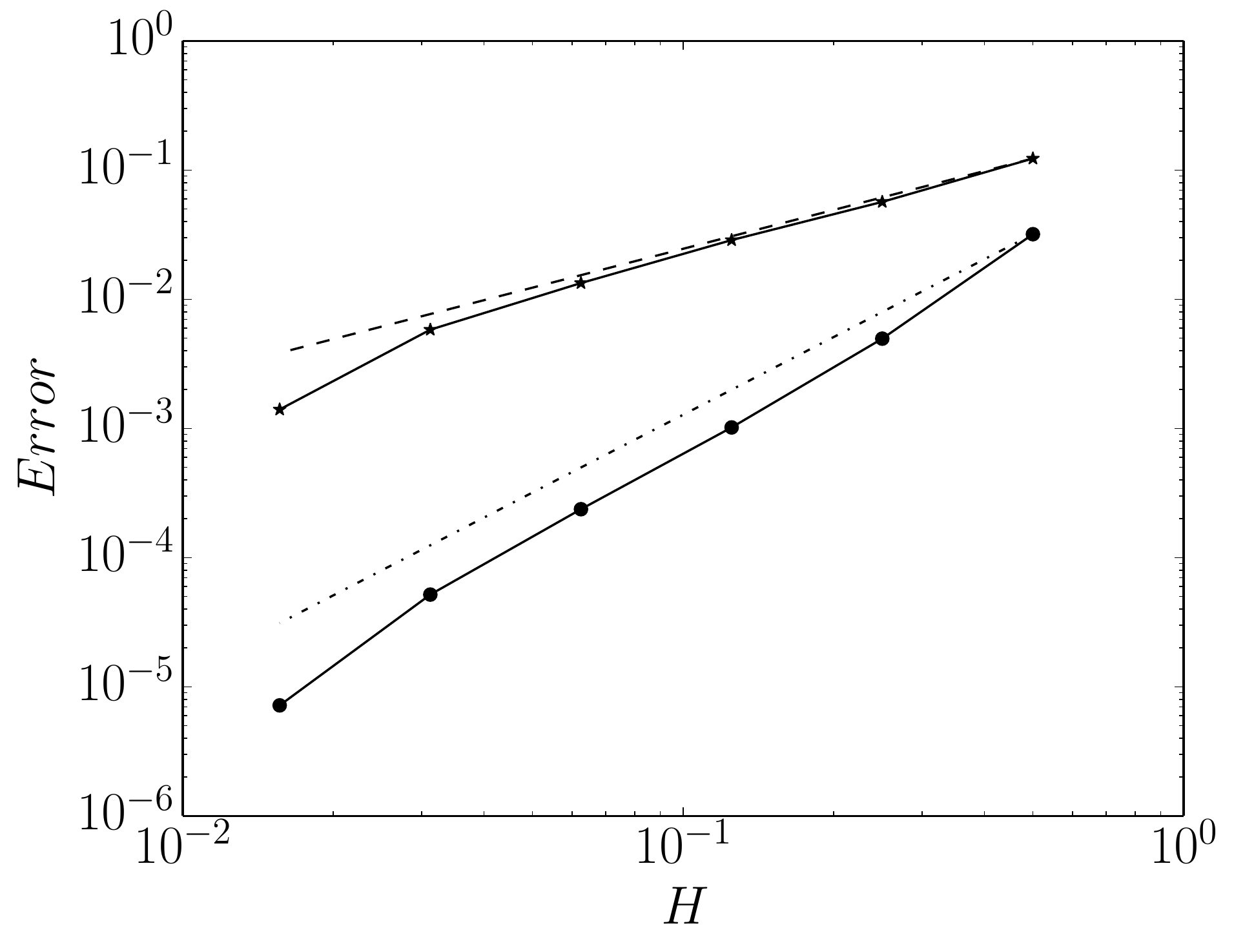}
    \includegraphics[width=.43\linewidth]{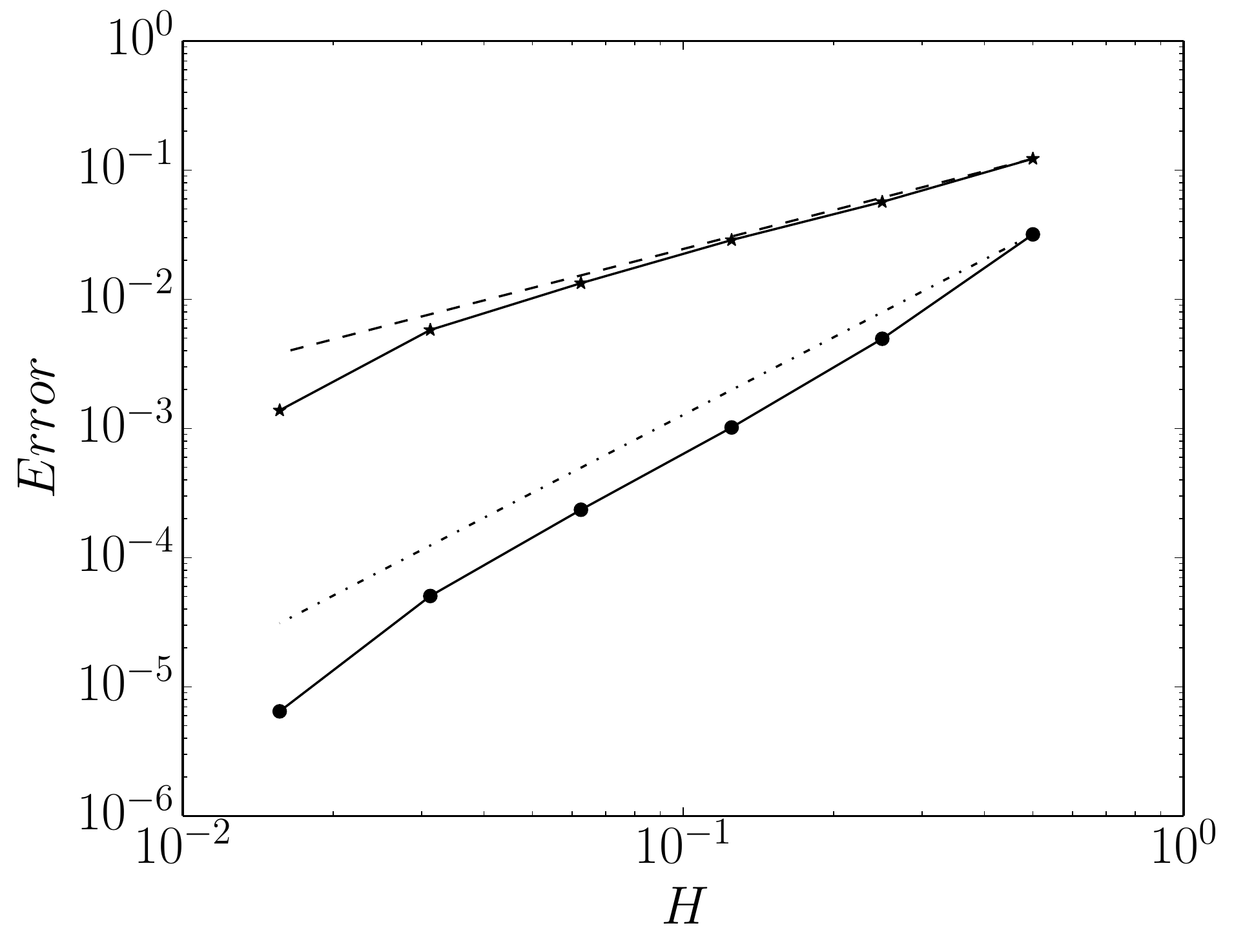}
    \includegraphics[width=.43\linewidth]{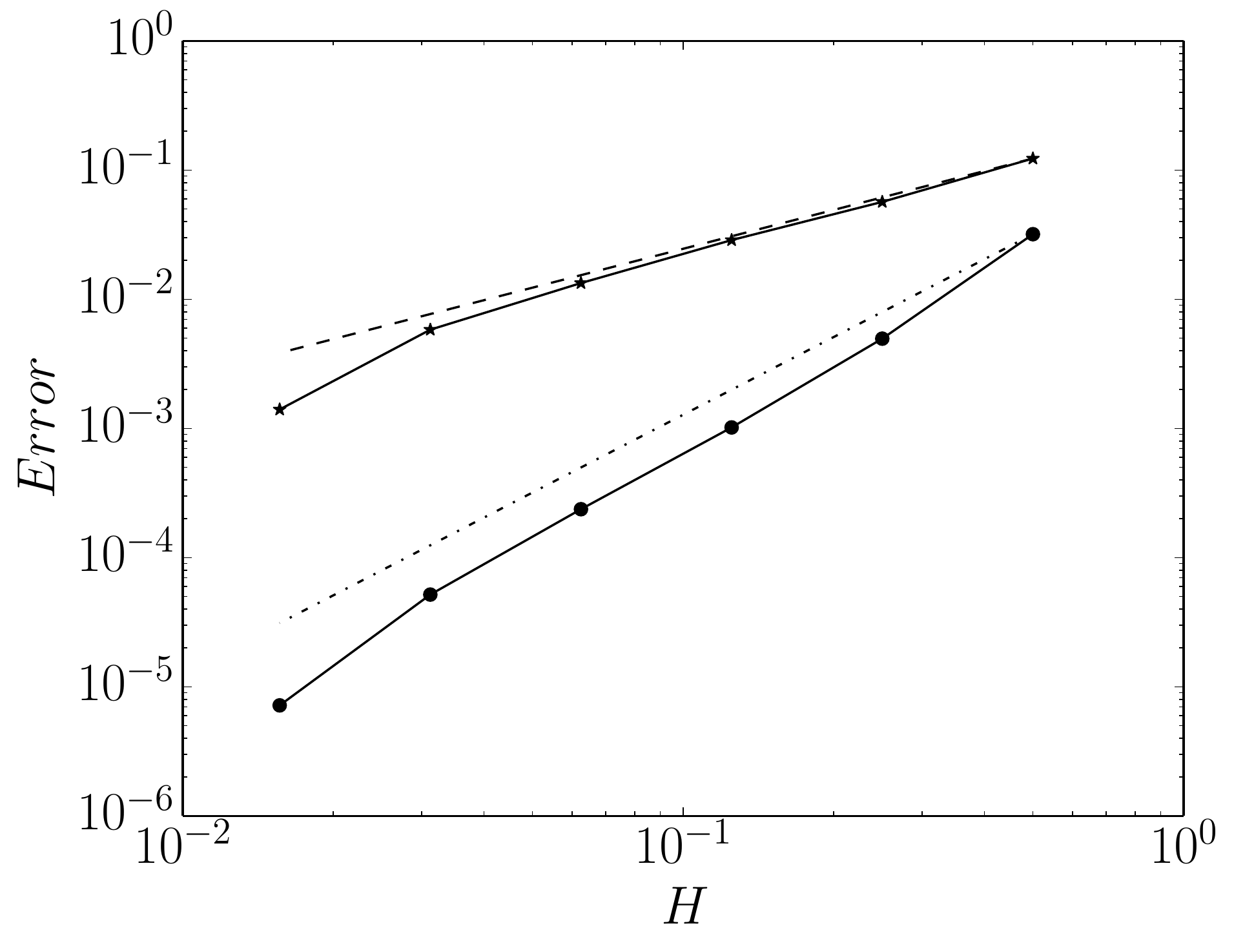}
    \includegraphics[width=.8\linewidth]{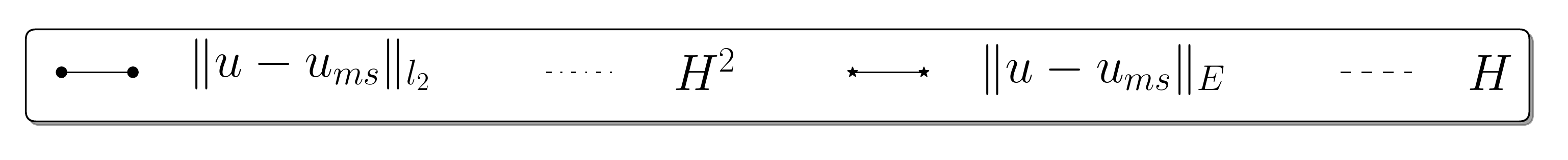}
    \caption{Residual of the LOD solver on a perturbed network with homogeneous (left) and random (right) coefficients. The upper graphs are for the fixed boundary problem and for the lower graphs displaced boundary problem. Optimal convergence rates are obtained in all of the experiments.}
    \label{fig:perturbed_test}
\end{figure}

The first setup is similar to one already confirmed \cite{DLOD}, where a perturbed network with uniform and random coefficients is used as the network model. This perturbed network is composed of $(2^{9}+1)^2$ network nodes. The results of testing this setup are presented by plotting the error of the multiscale approximation against the grid size $H$ for both the fixed boundary problem and displacement problem in Figure \ref{fig:perturbed_test}.  The results of these experiments coincide with previous results \cite{DLOD}, where the optimal convergence is observed for both perturbed networks with uniform and random coefficients. The randomization of force coefficients does not seem to affect the convergence rate in either problem, suggesting that coefficient variation does not affect the multiscale method's performance for discrete network models. 
\subsection{Structured and perturbed networks with paper-like coefficients}
\begin{figure}[H]
    \centering
    \includegraphics[width=.43\linewidth]{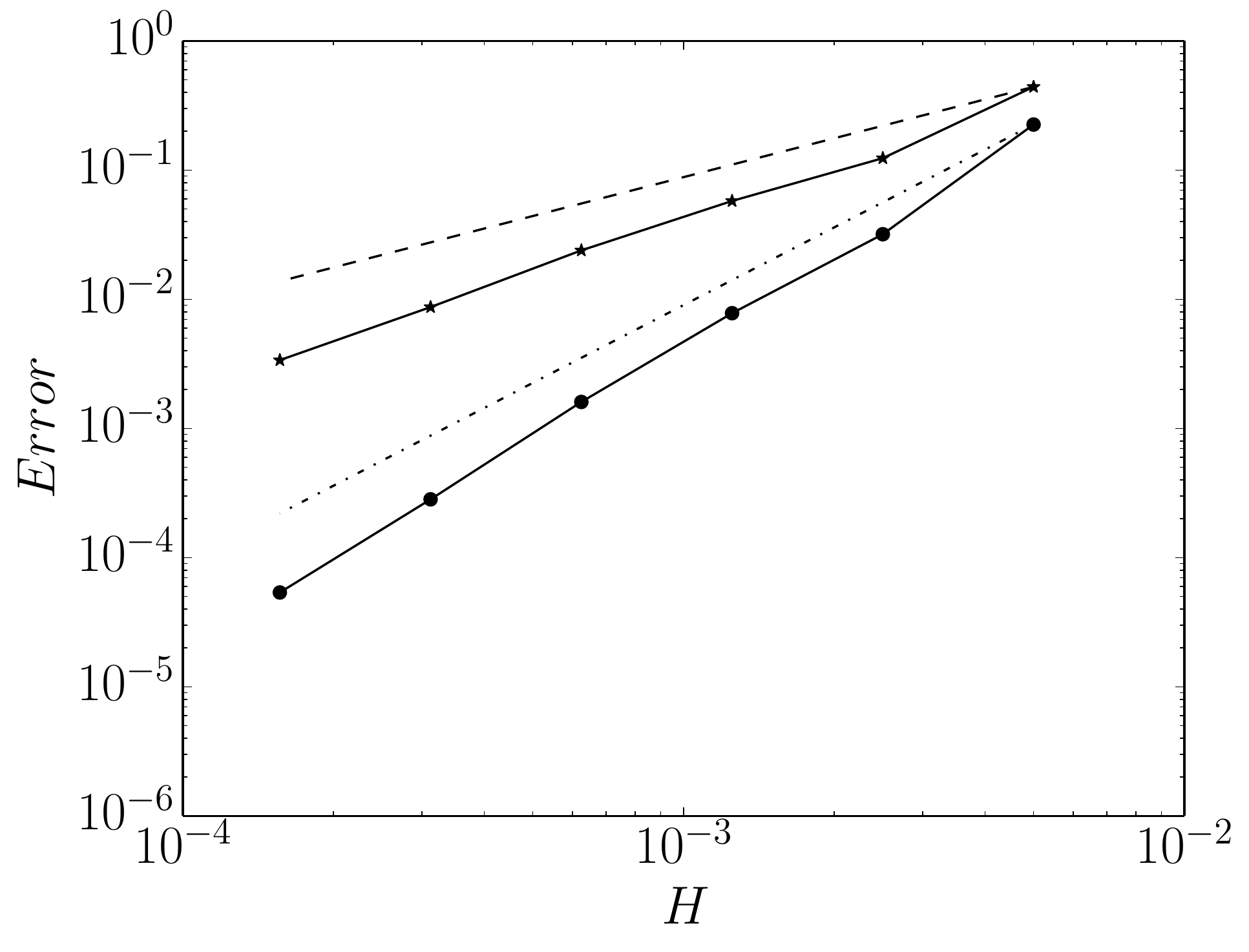}
    \includegraphics[width=.43\linewidth]{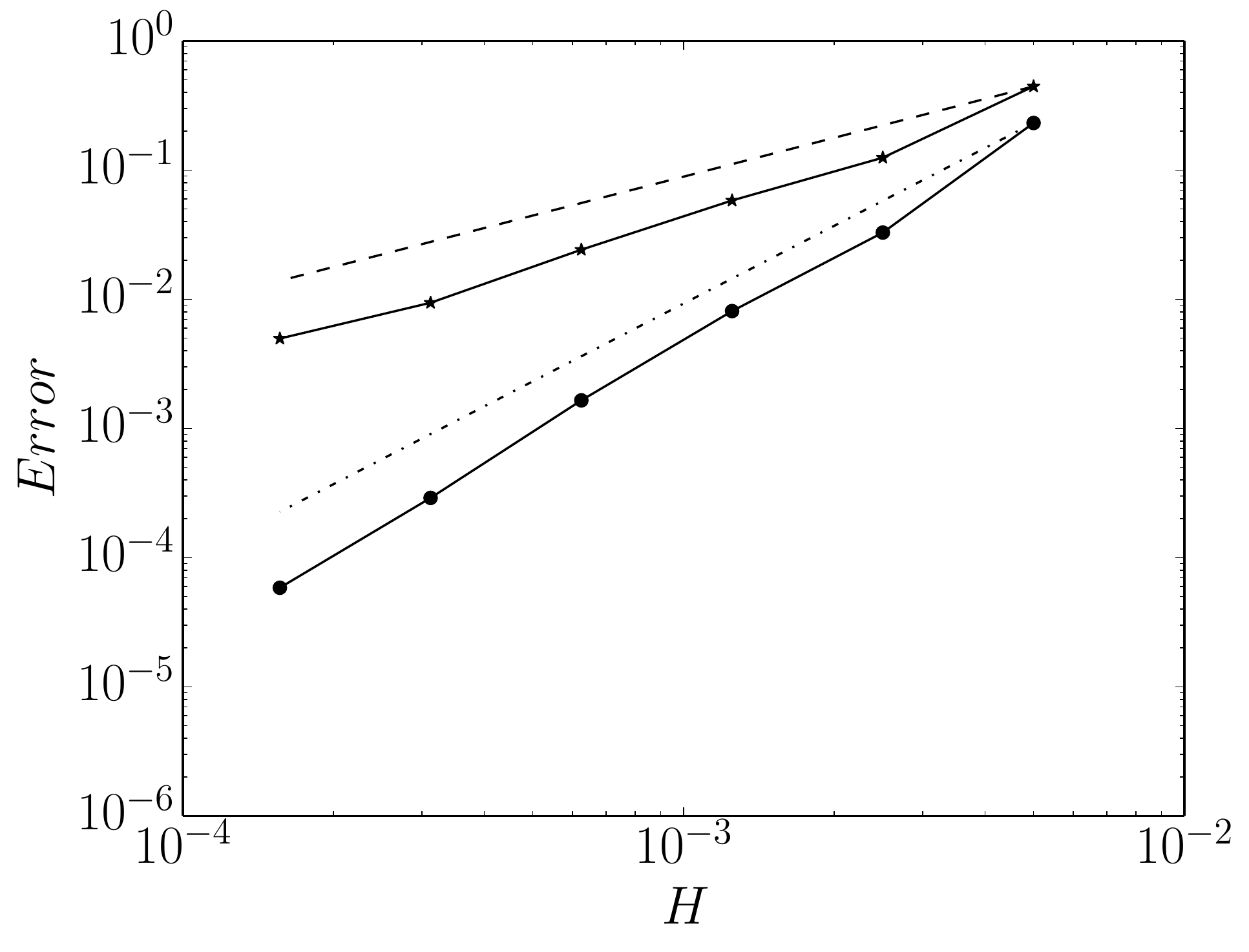}
    \includegraphics[width=.43\linewidth]{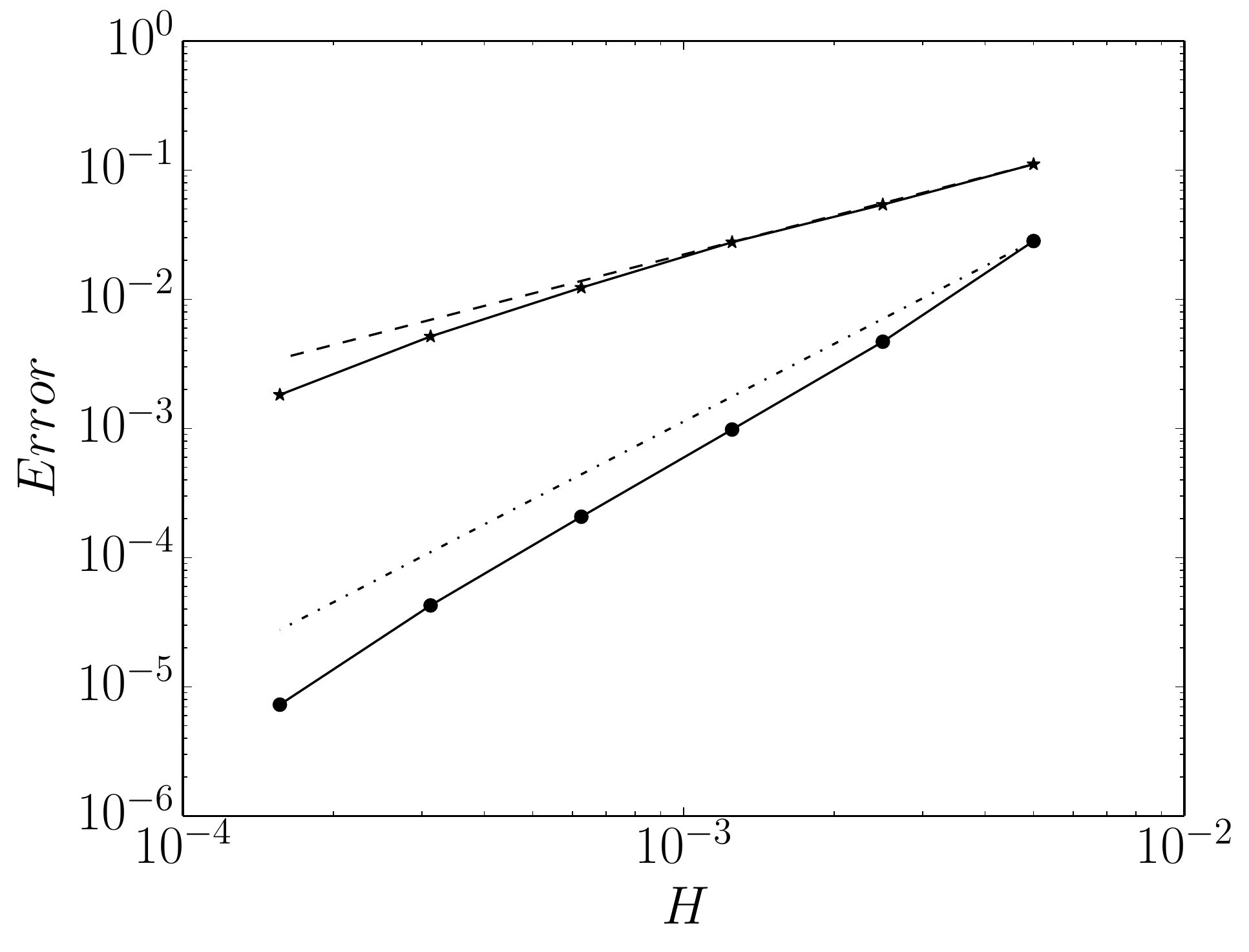}
    \includegraphics[width=.43\linewidth]{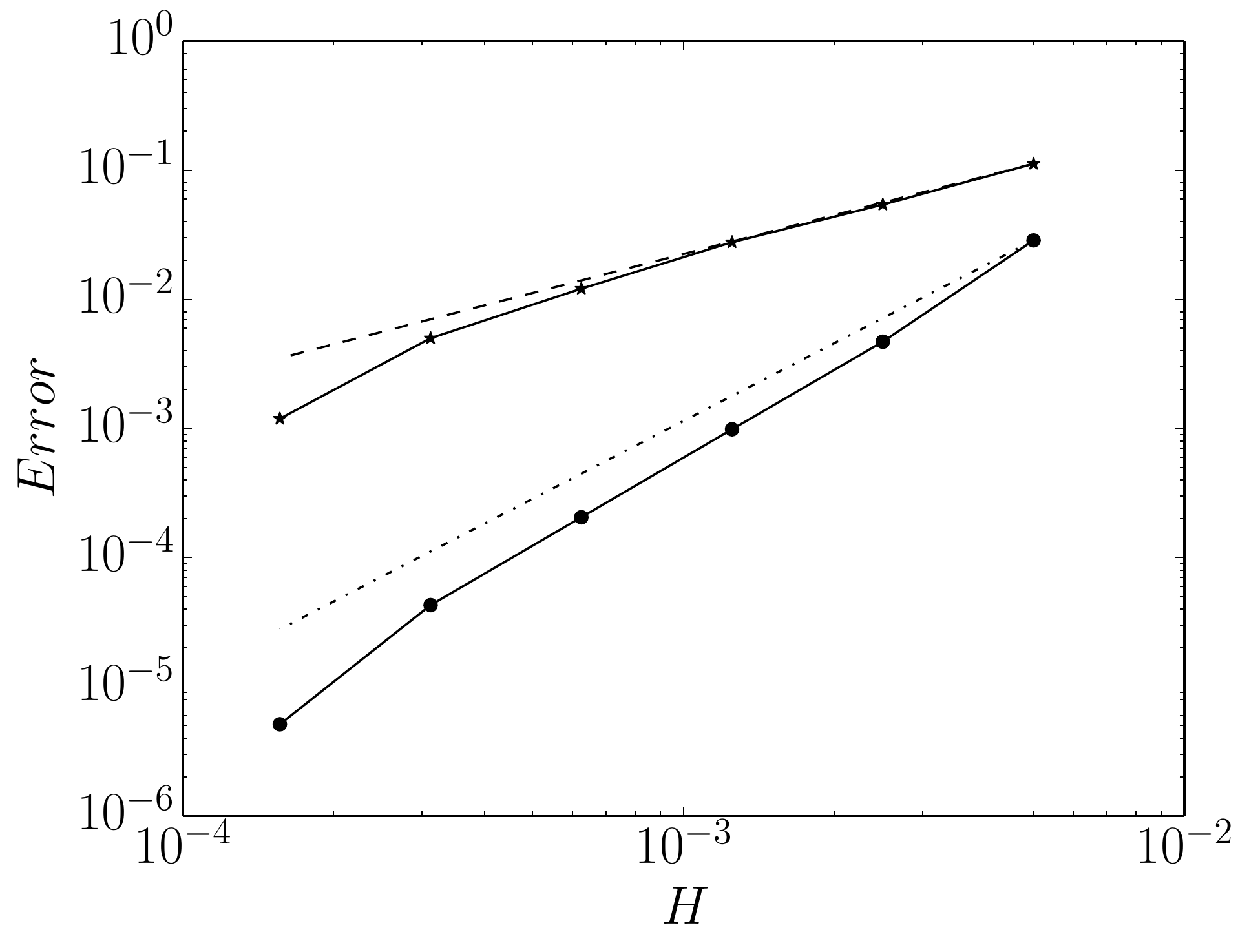}
    \includegraphics[width=.8\linewidth]{Legend.pdf}
    \caption{Residual of the LOD solver used on the displaced boundary problem. The left graphs show the convergence rates for structured networks and the right graphs the convergence rates for perturbed networks. The upper graphs are for the fixed boundary problem and the lower the displaced boundary problems. In these tests, the network models have stiff homogeneous, fiber, coefficients. The results show optimal convergence rates.}
    \label{fig:perturbedfiber_tests}
\end{figure}

The second setup is made to validate the method for stiff network models. In this setup, a structural and a perturbed network with $(2^{9}+1)^2$ network nodes are generated on the square domain $[0,0.01]^2$ with coefficients similar in scale to real paper. The change of domain is to make the edge density more similar to that of real paper. The results of the numerical experiments on these models are presented in Figure \ref{fig:perturbedfiber_tests} for both the fixed boundary problem and the displaced boundary problem. The convergence results for the perturbed network are similar to the results in the previous section. For the structured network, the convergence rates vary less and follow a more clear convergence path. Optimal convergence rates are observed in these tests.   
\subsection{Unordered fiber networks with paper-like coefficients}
\begin{figure}[H]
    \centering
    \includegraphics[width=.43\linewidth]{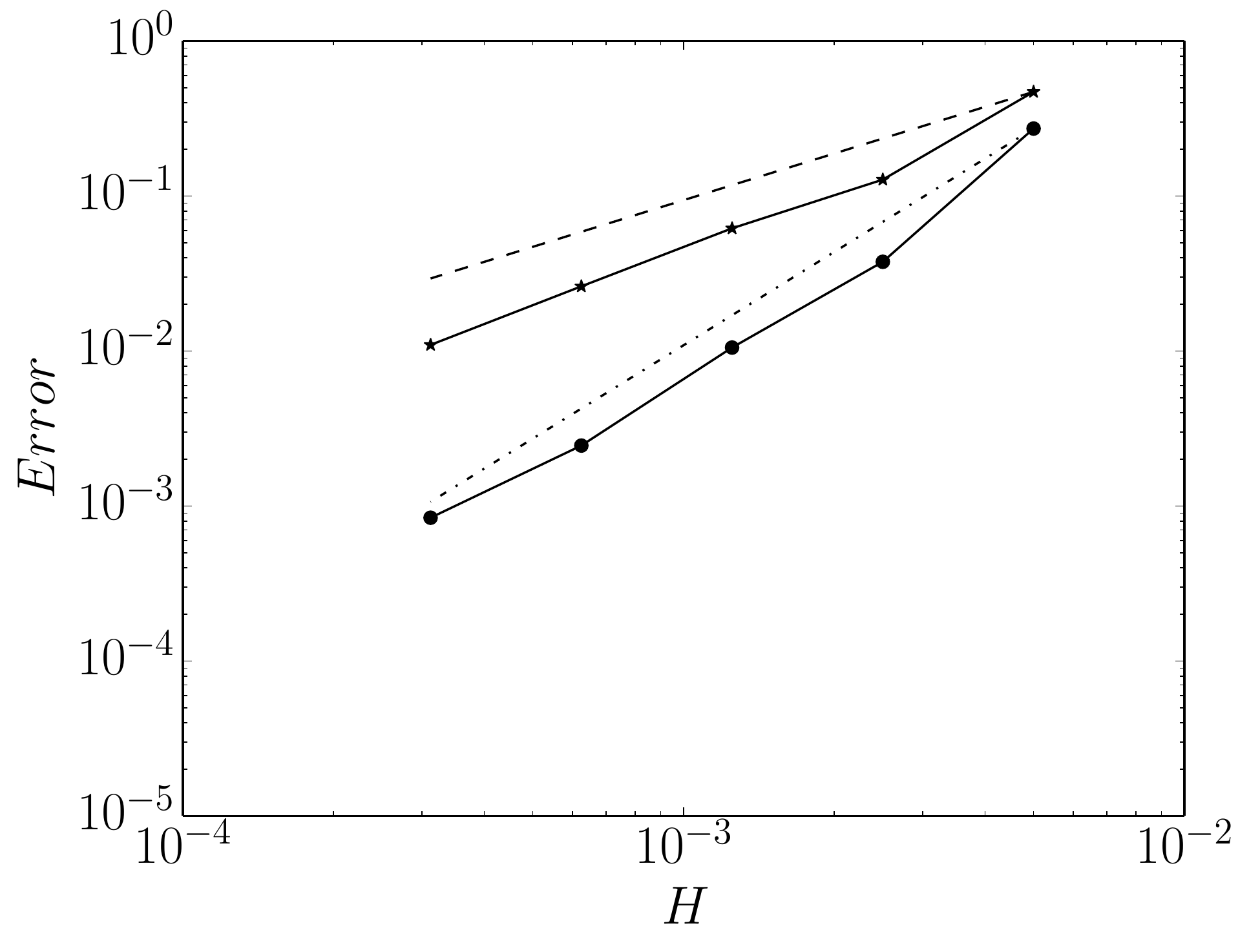}
    \includegraphics[width=.43\linewidth]{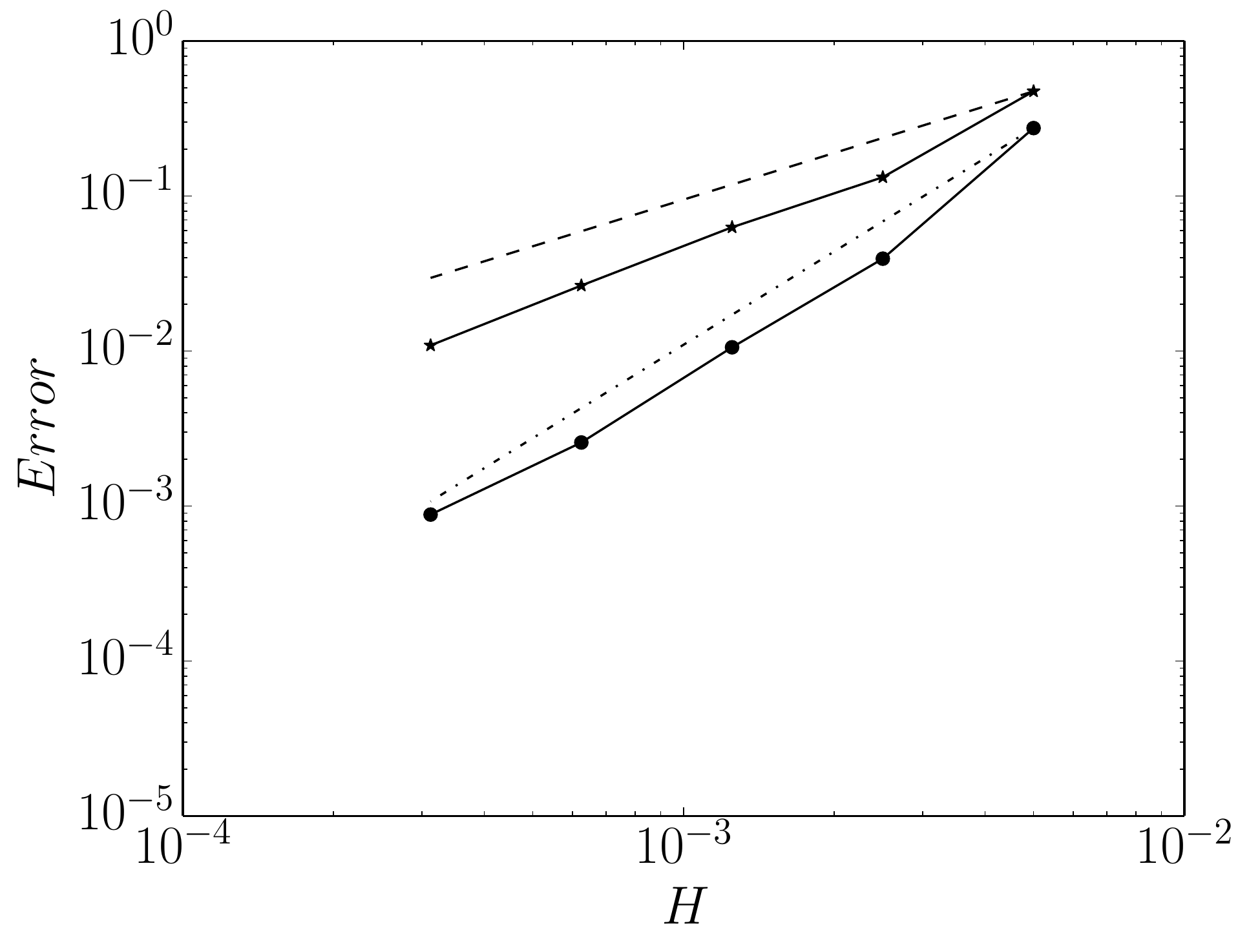}
    \includegraphics[width=.43\linewidth]{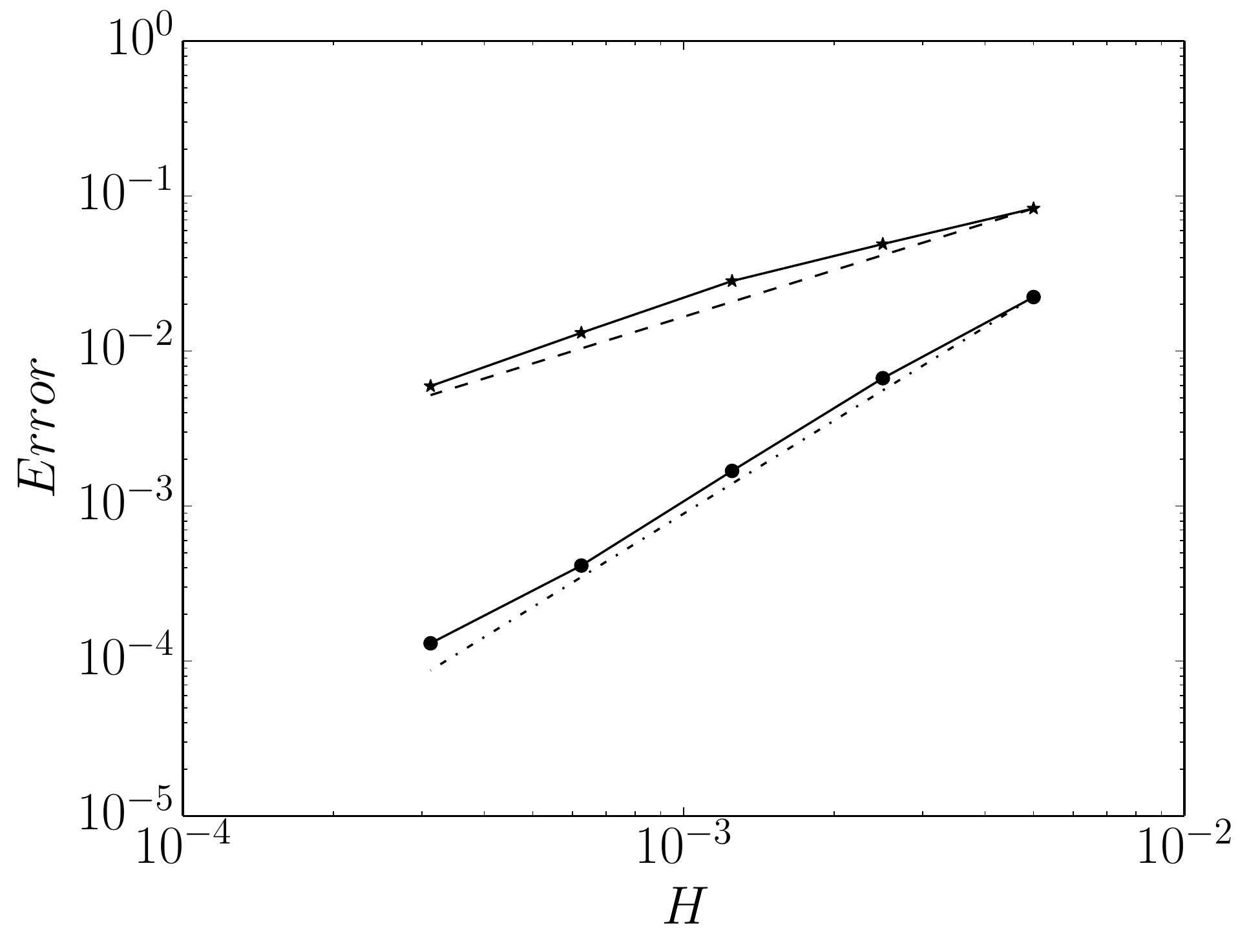}
    \includegraphics[width=.43\linewidth]{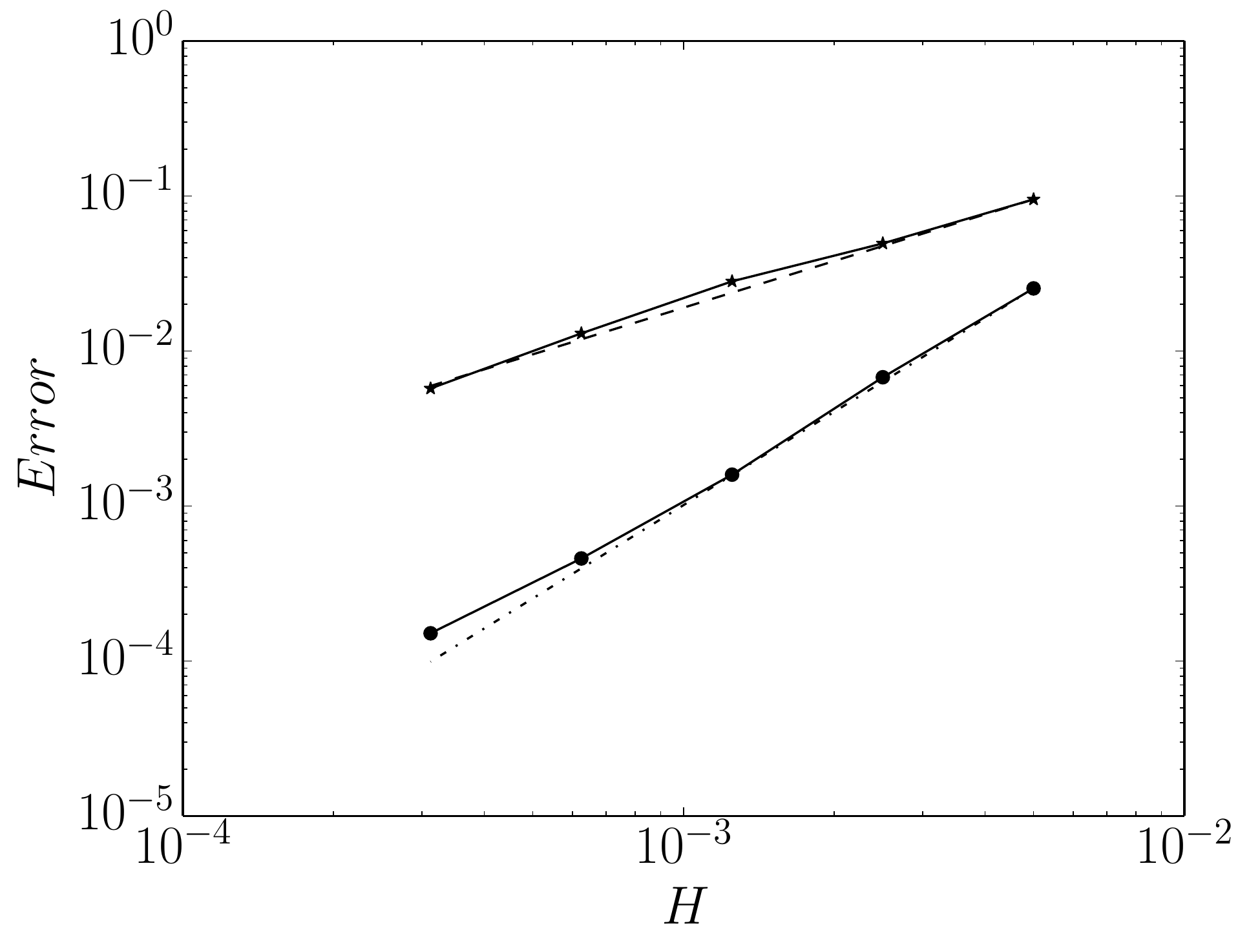}
    \includegraphics[width=.80\linewidth]{Legend.pdf}
    \caption{Errors of the LOD method's approximations on an unordered fiber network. Left graphs are with homogeneous fiber coefficients and the right graphs are with randomized fiber coefficients. The two top graphs are for the fixed boundary problem and the two lower graphs are for the displaced boundary problem. In the four experiments, optimal convergence results can be observed.}
    \label{fig:birdsnest_tests}
\end{figure}

In the last setup, a network model with similar geometrical and physical properties to real paper is analyzed. The model is an unordered fiber network with fiber coefficients on the domain $[0,0.01]^2$ with both uniform and random fiber coefficients. This model is comparable to a square pieces of paper with $12g/m^2$ grammage. This grammage on the given domain results in about two hundred thousand nodes which is similar to the two other setups. Because of the unstructured nature of these network models, the finest coarse grid analyzed is two times coarser than in the previous two setups. The convergence results for the fixed boundary and displacement problem for this model using the LOD method are presented in Figure \ref{fig:birdsnest_tests}. In the experiments conducted, optimal convergence rates are obtained, as in the previous results.
\section{CONCLUSIONS AND FUTURE WORK}
The LOD method has been numerically confirmed to handle network models with highly oscillating coefficients in two dimensions. Numerical experiments show that optimal convergence results are obtainable for network models using a high degree of localization. The models tested ranged from structured, semi-structured to unstructured networks. Structured and semi-structured, perturbed, networks have been numerically confirmed before with random parameters \cite{DLOD} and in this work the validation is extended to include stiff homogeneous and stiff random parameters. Unstructured stiff networks is a new family of networks numerically confirmed to work optimally with the LOD method. The network model used for evaluating the unordered network is not trivial, as it is inspired by a sheet of paper with edges representing segments of individual paper fibers. The LOD method can handle these highly irregular models with both stiff homogeneous and stiff random coefficients, optimally, with a high degree of localization. 

The goal of this ongoing project is to apply the method on verified paper models. With such a method and model, tensile strength and other physical properties can be obtained through simulation. For this to be possible, the method needs to work in three dimensions, which is the next step in the numerical verification process. Once three-dimensional network models can be handled, it will be possible to work with models generated in paper fabrication simulations or from experimental data. At this point, it seems like the local decomposition method is a good candidate for these types of models, but the conclusion is based solely on numerical results. Mathematical theory needs to be developed for this application of the LOD method before any real conclusions can be made. This theoretical development might reveal potential restrictions, but also inspire a wider range of applications.
\subsection*{Acknowledgements}
This work is a part of the ISOP (Innovative Simulation of Paper) project which is performed by a consortium consisting of Albany International, Stora Enso and Fraunhofer-Chalmers Centre. The main author is partially funded by the Swedish Foundation for Strategic Research (SSF).

\end{document}